\documentstyle[12pt]{amsart}
\newtheorem{prop}{Proposition}
\newtheorem{th}[prop]{Theorem}

\newtheorem{lem}[prop]{Lemma}

\theoremstyle{definition}
\newtheorem{ack}{Acknowledgments} 
\theoremstyle{remark}
\newtheorem{rem}{Remark} 
\newcommand{\bbN}{{\Bbb{N}}}
\newcommand{\bbR}{{\Bbb{R}}}
\newcommand{\bbK}{{\Bbb{K}}}
\newcommand{\calB}{{\cal{B}}}
\newcommand{\calD}{{\cal{D}}}
\newcommand{\calF}{{\cal{F}}}

\newcommand{\m}{\mu}
\newcommand{\n}{\nu}

\newcommand{\s}{\sigma}
\newcommand{\f}{\varphi}
\newcommand{\h}{\chi}

\newcommand{\om}{\omega}
\newcommand{\Om}{\Omega}

\newcommand{\supp}{\operatorname{supp}}
\newcommand{\card}{\operatorname{card}}

\newcommand{\sgn}{\operatorname{sgn}}

\newcommand{\dint}{\displaystyle\int}
\newcommand{\disp}{\displaystyle}
\newcommand{\lb}{\label}

\newcommand{\emp}{\emptyset}

\newcommand{\lra}{\longrightarrow}

\newcommand{\wtw}{if and only if }
\newcommand{\buo}{without loss of generality }

\makeatletter
\def\@currentlabel{2.1}\label{e:dispaa}
\def\@currentlabel{2.21}\label{e:dispau}
\def\@currentlabel{2.22}\label{e:dispav}
 \def\@currentlabel{2.23}\label{e:dispaw}
\def\@currentlabel{2.24}\label{e:dispax}
\makeatother

\makeatletter
\def\alphenumi{%
  \def\theenumi{\alph{enumi}}%
  \def\p@enumi{\theenumi}%
  \def\labelenumi{(\@alph\c@enumi)}}
\makeatother

\begin{document}

\title{Isometries of Hilbert space valued  function
spaces}
\author{Beata Randrianantoanina}
\address{Department of Mathematics \\ The University  of Texas at  
Austin
\\
         Austin, TX 78712}
\email{brandri@@math.utexas.edu}
\maketitle

\begin{abstract} Let  $X$ be a (real or complex)
rearrangement-in\-va\-riant  function space on $\Om$ (where $\Om =  
[0,1]$ or $\Om \subseteq \bbN$) whose norm is not proportional to  
the
$L_2$-norm. Let $H$ be a separable Hilbert space.
We characterize surjective isometries of $X(H).$ We prove that if $T$  
is such an isometry then
there exist Borel maps $a:\Om\to\bbK$ and  $\sigma:\Om\lra\Om$ and a  
strongly measurable operator map $S$ of $\Om$ into $\calB(H)$
so that for almost all $\om$ $S(\om)$ is a surjective isometry of  
$H$   

and for any $f\in X(H)$
$$Tf(\om)=a(\om)S(\om)(f(\sigma(\om))) \text{  a.e.}$$ 

As a consequence we obtain a new proof of characterization of  
surjective isometries in complex rearrangement-invariant function
spaces.
\end{abstract}

\section{Introduction}

We study isometries of Hilbert space valued rearrangement-invariant  
function
spaces $X(H)$, where $\dim H\ge 2$ and $H$ is separable..
Our results are valid for both symmetric sequence spaces
and nonatomic rearrangement-invariant function spaces on $[0,1]$ with  
norm
not proportional to $L_2$
but they are new only in the nonatomic case.
If $X$ is a sequence space, even not necessarily symmetric,
Theorem~\ref{isoxl2} is a special case of a much more general result  
of
Rosenthal \cite{R86} about isometries of Functional Hilbertian Sums.
We include here the case of $X$ being a symmetric sequence space  
since the
proof is essentially the same as when $X$ is a nonatomic  
rearrangement-invariant
function space, and also our techniques are much simpler than those
developed in \cite{R86}.

 Spaces of the form $X(H)$ appear naturally in the theory of Banach  
spaces
 (see \cite[Chapter~2.d]{LT2}). In particular, if $X$
 is rearrange\-ment-invariant (with Boyd indeces $1<p_X \le q_X  
<\infty$) then $X(L_2)$ is isomorphic to $X$  
(\cite[Proposition~2.d.4]{LT2}) and this plays an important role in  
the study of uniqueness of unconditional bases in $X.$

Isometries of Hilbert space-valued function spaces have been studied  
by many authors. Cambern \cite{C74} (1974) characterized isometries  
of $L_p(L_2)$ in the complex case (see also an alternative proof of  
Fleming and Jamison \cite{FJ77}). Isometries  of $L_p(L_2)$ in both  
real and complex cases are described (among other spaces) in the  
general paper of Greim \cite{Greim} (1983).  In 1981 Cambern  
\cite{C81} described isometries of, both real and complex,  
$L_{\infty}(L_2).$ In 1986 Jamison and Loomis
\cite{JL}
gave the characterization of isometries in complex Hilbert  
space-valued nonatomic Orlicz spaces $X(L_2).$ Also there have been  
done a number of studies of various $L_2$-valued analytic function  
spaces. For a fuller discussion of literature we refer the reader to  
the forthcoming survey of Fleming and Jamison \cite{FJs}.

We use a method of proof which is designed for spaces over $\bbR$,  
but clearly
complex linear operators $T:X(H) \to X(H)$ can be always considered  
as real linear
operators acting on $X(H)(\ell_2^2)$ and therefore our results are  
valid also
in complex case.

Moreover Theorem~\ref{isoxl2} with $H=\ell_2^2$ may be viewed as a  
statement
about the form of isometries of complex rearrangement-invariant  
spaces. Thus we give a new proof of the fact that all surjective  
isometries on $X$ can be represented as weighted composition  
operators i.e. if $T$ is such an isometry then there are Borel maps  
$a, \sigma$ such that $Tf = a f\circ \sigma$ for all $f$ in $X$
(cf. \cite{Z77}, \cite{Z80} for nonatomic spaces and \cite{T} for
sequence spaces).

\section{Preliminaries} \lb{sec2}

We follow standard notations as in \cite{LT2}.

In the following $H$ denotes a separable Hilbert space with $\dim H  
\ge 2.$ If we want to stress that we restrict our attention to the  
case when $\dim H = \infty$ we will write $H=\ell_2.$

If $X$ is a K\"othe function space (\cite[Definition 1.b.17]{LT2}) we  
denote  by  $X'$ the {\bf K\"othe dual} of $X;$   thus $X'$ is the
K\"othe space of all $g$ such that $\int|f||g|\,d\mu<\infty$ for  
every $f\in X$ equipped with the
norm $\|g\|_{X'}=  \sup_{\|f\|_X\le 1}\int |f||g|\,d\mu.$
Then $X'$ can be regarded as a closed subspace of the dual
$X^*$ of $X$.

If $X$ is a K\"othe function space on $(\Om_1,\mu_1)$ and $H$ is a
separable Hilbert space on $(\Om_2,\m_2),$ we will denote by
$\bold{X(H)}$ \lb{defxy} the K\"othe function space on
$(\Om_1 \times \Om_2,\mu_1 \times \m_2)$ with a following norm
$$\|f(\om_1,\om_2)\|_{X(H)} = \| \ \|f(\om_1,\cdot)\|_2 \ \|_X.$$

This definition coincides with the notion of $H$-valued Bochner  
spaces.

 It is well-known that
$(X(H))^* = X^*(H),$ and the space $(X(H))' \subset X^*(H)$ can be
identified with the space of functions $\f:\Om_1 \lra H$ such that  
for
every $y \in H$ the map $\om_1 \longmapsto \langle\f(\om_1),y\rangle$  
is
measurable and the map $\f_\#:\om_1 \longmapsto \|\f(\om_1)\|_{H}$  
belongs
to $X'.$ The operation of $\f$ on $X(H)$ is given by
$$\f(f) = \int_{\Om_1} \langle\f(\om_1),f(\om_1)\rangle\  
d\m_1(\om_1)$$
for any $f \in X(H).$  Thus $(X(H))' = X'(H).$

For any function $f\in X(H)$ we define the map $f_\# :  
\omega_1\to\bbR$
by $f_\#(\omega) = \|f(\omega)\|_{H}$.
Then $f_\# \in X$. We say that functions $f,g \in X(H)$ are
{\bf disjoint in a vector sense} if $f_\#$ and $g_\#$ are disjointly
supported, i.e., $f_\# (\omega) \cdot g_\# (\omega)=0$ for a.e.  
$\omega\in
\Omega_1$. We say that an operator $T:X(H)\to X(H)$ is {\bf  
disjointness
preserving in a vector sense} if $(Tf)_\# \cdot (Tg)_\# =0$ whenever  
$f_\#
\cdot g_\#=0$.

We will say that an operator $T:X(H) \to X(H)$ has {\bf a canonical  
vector
form} if there exists a nonvanishing Borel function $a $ on $\Om$  
(where $\Om
= [0,1]$ if $X$ is nonatomic or $\Om \subset \bbN$ if $X$ is a
sequence space) and an invertible Borel map $\sigma:\Om\lra\Om$ such
that, for any Borel set
$B \subset \Om,$ we have
$\m(\sigma^{-1}B)=0$ if and only if $\m(B)=0$ and
a strongly measurable map $S$ of $\Om$ into $\calB(H)$ (i.e. for each  
$h\in H$
the mapping $\om \mapsto S(\om)h$ is measurable)
so that $S(t)$ is an isometry of $H$ onto itself for almost all $t$
and
$$Tf(t)=a(t)S(t)(f(\sigma(t))) \text{  a.e.}$$ for any $f\in X(H).$

Note that the name ``a canonical vector form'' is introduced here  
only for the
purpose of this paper --- we do not know the standard name for this  
type
of operator.  We will need the following simple observation
(cf. \cite[ Lemma 2.4]{KR})

\begin{lem}\label{2.4'}
Suppose that $T:X(H)\to X(H)$ is an invertible operator which has a  
canonical
vector form. Then $T':X'(H)\to X'(H)$ exists and has a canonical  
vector form.
\end{lem}

\begin{pf}
Operator $T$ has an representation
$$Tf(\omega_1) = a(\omega_1) S(\omega_1)(f(\sigma(\omega_1)))$$
where $a,S,\sigma$ satisfy the above conditions for canonical form  
and
moreover  $a$ is nonvanishing and $\sigma$ is
an invertible Borel map with $\mu(\sigma^{-1}B)=0$ if and
only if $\mu(B)=0.$ Let $v$ be the Radon-Nikodym derivative
of the $\sigma$-finite measure $\nu(B)=\mu(\sigma^{-1} B).$

Then for $f\in X(H),g\in X'(H)$ we have
\begin{align*}
g(Tf) &=  \int_{\Omega_1}\langle g(\omega_1),Tf(\omega_1)\rangle  
\,d\mu (\omega_1)=\\
&= \int_{\Omega_1}\langle g(\omega_1),a(\omega_1)S(\omega_1)
(f(\sigma(\omega_1)))\rangle\,d\mu (\omega_1)=\\
&= \int_{\Omega_1}\langle a(\omega_1)(S(\omega_1))'(g(\omega_1)),
f(\sigma(\omega_1)) \rangle\,d\mu(\omega_1)=\\
&= \int_{\Omega_1}\langle  
a(\sigma^{-1}(\omega_1))(S(\sigma^{-1}(\omega_1)))'
(g(\sigma^{-1}(\omega_1))),f(\omega_1)\rangle v(\omega_1)\,  
d\mu(\omega_1)\ ,
\end{align*}
since $(S(\om))^* = (S(\om))'.$

Thus $T^*g \in X'(H)$ and
$$T'g(\omega_1) = a(\sigma^{-1}(\omega_1))v(\omega_1)
(S(\sigma^{-1}(\omega_1)))' g (\sigma^{-1}(\omega_1))\ \text{  
a.e.}$$
Clearly the map $\omega_1\mapsto S(\sigma^{-1}(\omega_1))'$ is  
strongly
measurable and thus $T'$ has a canonical vector form.
\end{pf}

A {\bf rearrangement-invariant function space} ({\bf r.i. space})
\cite[Definition 2.a.1]{LT2} is a K\"othe function space on  
$(\Om,\mu)$
which satisfies the conditions:
\begin{enumerate}
\item[(1)]  $X'$ is a norming subspace of $X^*.$
  \item[(2)] If $\tau:\Om\lra \Om$ is any
measure-preserving invertible Borel automorphism then $f\in
X$ if and only if $f\circ\tau\in X$ and
$\|f\|_X=\|f\circ\tau\|_X.$
\item[(3)] $\|\chi_{B}\|_X=1$ if $\mu(B) = 1.$
\end{enumerate}

Next we will quickly state a definition of Flinn elements.
For fuller description and proofs we refer to \cite{KR} and  
\cite{th}.

We say that an element $u$ of a K\"othe space $X$ is {\bf Flinn}  if
there exists an $f\in X^*$ such that $f\ne 0$ and for every $x\in X$  
and
$x^*\in X^*$ with $x^* (x) = \|x\|_X\cdot \|x^*\|_{X^*}$ we have
$$f(x)\cdot x^* (u)\ge 0\ .$$
We say that $(u,f)$ is a {\bf Flinn pair}.
We denote by ${\calF}(X)$ the set of Flinn elements in $X$.
We will need the following facts:

\begin{prop} \lb{fiso} (\cite[Proposition 3.2]{KR})

Suppose $U:X\lra Y$ is a surjective isometry.
Then  $U(\calF (X)) = \calF (Y)$; furthermore if \/ $(u,f)$
is a Flinn pair then \/ $(U(u),(U^*)^{-1}f)$ is a Flinn pair.
\end{prop}

\begin{th}   \lb{flinn} (Flinn, \cite[Theorem 1.1]{R84},  
\cite[Theorem 3.3]{KR})

Let $X$ be a Banach space and $\pi$ be a
contractive projection on $X$ with range $Y.$  Suppose $(u,f)$ is a  
Flinn
pair in $X$.  Suppose $f\notin Y^{\perp}.$  Then $\pi(u)\in  
\calF(Y).$
\end{th}

\begin{th}  \lb{kw} (\cite[Theorem 4.3]{KR})

Suppose $\mu$ is nonatomic and suppose
$X$ is an order-continuous  K\"othe function space on  
$(\Omega,\mu).$
Then $u\in X$ is a Flinn element
 if and only if there is a nonnegative
function $w\in L_0(\mu)$ with supp $w=$ supp $u=B,$ so that:\newline
$(a)$ If $x\in X(B)$
then
${\disp \|x\|=(\int |x|^2w\,d\mu)^{1/2}.}$ \newline
and\newline
$(b)$ If $v\in X(\Omega\setminus B)$ and $x,y\in X(B)$ satisfy
$\|x\|=\|y\|$ then
$\|v+x\|=\|v+y\|.$
\end{th}

The last fact about Flinn elements that we will
need is a reformulation of Calvert and Fitzpatrick's
characterization of $\ell_p-$spaces \cite{CF86}:

\begin{th} \lb{cf}
Suppose that $X$ is a sequence space with $\dim X = d <\infty, d\ge  
3,$ and basis $\{e_i\}_{i=1}^d.$ Suppose that every element $u$ of  
$X$ with support on at most two coordinates is Flinn in $X,$ i.e.
$$\{u\in X\ : \ u = a_i e_i + a_je_j \text{ for some $i,j\le d, \ \  
a_i,a_j \in\bbR$}\} \subset \calF(X).$$
Then $X= \ell_p^d$ for some $1\le p \le \infty.$
\end{th}

\begin{pf}
By Lemma~1.4 of \cite{R86} $(u,f)$ is a Flinn pair in $X$ \wtw the
projection $P$ defined by $P(x) = x - f(x)u$ has norm 1 in $X.$  
Hence, if $(u,f)$ is a Flinn pair in $X$ then there is a projection  
of norm 1 onto the hyperplane $\ker f\subset X.$

It is also clear from the definition that if $(u,f)$ is a Flinn pair  
in $X$ then $(f,u)$ is a Flinn pair in $X'.$ Therefore there exists a  
projection of norm 1 onto $\ker u\subset X'$ for every $u$ with  
support on at most two coordinates.
But then Theorem~1 of \cite{CF86} asserts that if
$d\ge 3$ then $X'= \ell_q^d$ for some $1\le q \le \infty.$
Thus $X= \ell_p^d.$
\end{pf}

Finally let us introduce the following notation.

Suppose that $X$ is a nonatomic r.i. \ space on $[0,1]$ and $n$ is a  
natural
number.   Let
$e^n_i=\chi_{((i-1)2^{-n},i2^{-n}]}$ for $1\le i\le 2^n.$
Denote $X_n=[e^n_i:1\le i\le 2^n].$ If $\dim X < \infty$ then, for  
the uniformity
of notation, we will use $X_n = X$ for any $n \in\bbN.$
 Notice that $X_n^*$ can be identified
naturally with $X'_n.$

We now need to introduce a technical definition.  We will
say that an r.i. space $X$ has {\bf property} ${\bold(P)}$ if for  
every
$t>0,$
\begin{enumerate}
\item[a)] $\| \chi_{[0,\frac12]}\|_X < \| \chi_{[0,\frac12]}+
t\chi_{[\frac12,1]} \|_X $ if $X$ is a nonatomic function space on  
$[0,1];$ or
\item[b)]
$ \|e_1\|_X < \|e_1+te_2\|_X$ if $X$ is a sequence space with basis
$\{e_i\}_{i=1}^{\dim X}.$
\end{enumerate}
We say that
$X$ has property ${\bold(P')}$ if $X'$ has property (P).

Notice that, clearly, if $X$ has property $(P)$ (resp. $(P')$) then  
for every
$n\in \bbN$ \/ $X_n$  has property $(P)$ (resp. $(P')$).

\begin{lem} \lb{4.2} (\cite[Lemma 5.2]{KR})

Any r.i. space $X$ has at least one of
the properties $(P)$ or $(P').$ \end{lem}

The reason for introducing property $(P)$ is the following fact which  
will be
important for our applications.

If 
$v\in X_n(H)$ then $v=(v_i)_{i=1}^{2^n},$ where $v_i\in H$ for all  
$i$ and 
$v_i = (v_{i,j})_{j=1}^{\dim H}.$ 
Similarly for 
$f\in X_n'(H),$ \/ $f=(f_i)_{i=1}^{2^n},$ and   
$f_i = (f_{i,j})_{j=1}^{\dim H}\in H.$ In this notation we have:

\begin{lem}\lb{5.2}
Suppose that $X$ has property $(P')$ and $v\otimes f$ is a Flinn  
pair
in $X_n(H)$.
If $\|v_1\|_2 = |v_{11}|$ then $f_{11}\ne 0$.
\end{lem}

\begin{pf}

Assume that $f_{11} = 0$.
Then, since $v\otimes f\not\equiv 0$ there exists $i>1$ and $j\ge 1$  
such
that $f_{ij}\ne0$ and $v_{ij}\ne0$.
In fact $v_{ij} f_{ij}>0$ since $f(e_{ij})\cdot e_{ij}^*(v)\ge 0$.

Consider
$$e_{11}^* + te_{ij}^* \in X'_n (\ell_2^d)\ .$$
Then
$$\|e_{11}^* + te_{ij}^*\|_{X'_n(H)}
= \| e_1^* + te_i^*\|_{X'_n} > \|e_1^*\|$$
for all $t\ne 0$ since $X$ has $(P')$.
Hence for any $t\ne0$ if an element $(a_t e_{11} + b_t e_{ij})$ in
$X_n(H)$ is norming for $(e_{11}^* + te_{ij}^*)$ then $b_t\ne0$.
In fact $b_t\cdot t>0$.
Let us take $t= {-v_{11}\over 2v_{ij}}$.
Then $\sgn b_t = \sgn t = -\sgn (v_{11}\cdot v_{ij}) = -\sgn (v_{11}  
f_{ij})$.
Further:
\begin{align*}
&f(a_te_{11} + b_te_{ij})\cdot \left(e_{11}^* - {v_{11}\over 2v_{ij}}  
e_{ij}^*
\right) (v) = \\
&\qquad = b_t f_{ij} \cdot \left( v_{11} - {v_{11}\over 2v_{ij}}  
v_{ij}
\right) = \frac12 b_t\cdot f_{ij}\cdot v_{11} < 0\ .
\end{align*}
and the resulting contradiction with numerical positivity of  
$v\otimes f$
proves the lemma.
\end{pf}

\section{Main results} \lb{secl2}

We start with with an important (for us) proposition about the form  
of
Flinn elements in $X_n(H).$ In the case when $\dim H < \infty$ our
proof requires a certain technical restriction on the space $X,$
which is irrelevant in the case when $H = \ell_2.$ We present here
proofs for both cases since they are quite different. However, for
the application to Theorem~\ref{isoxl2} we need only to know the
validity of Proposition~\ref{5.2a}.

\begin{prop}\lb{5.2a}
Suppose that $X$ is an r.i.\ space  with property $(P')$, $\dim X \ge  
3$ and such that norm of $X$ is not
proportional to the $L_p-$norm for any $1\le p\le \infty$.
Then there exists $N\in\bbN,$ such that if $n\ge N$ and $u =  
(u_i)_{i=1}^{2^n} \in \calF(X_n(H))$ then there exists $1\le i_0\le  
2^n$ such that $\|u_i\|_2=0$ for all
$i\ne i_0$.
\end{prop}

\begin{rem}
Proposition~\ref{5.2a} can be also under\-stood as a state\-ment  
about the form of 1-co\-dimen\-sional hyperplanes in $X_n(H)$ which  
are ranges of a norm-1 projection.
\end{rem}

\begin{pf}
Let $n$ be big enough so that $X_n\ne \ell_p^{2^n}$, $1\le  
p\le\infty$.
Let $u\in {\calF}(X_n(H))$.
$$u= (u_i)_{i=1}^{2^n},\ \ \ \ \ \ \ \ \ \ \  u_i\in H.$$
Let $m= \card \{i:u_i\not\equiv 0\}$.
We want to prove that $m=1$.

By Proposition~\ref{fiso} we can assume without loss of generality  
that
$u_i\not\equiv 0$ for $i=1,\ldots,m$, $u_i=0$ for $i>m$ and $\alpha_1  
=
\|u_1\|_2 = \min \{\|u_i\|_2 : i=1,\ldots,m\}$.
Now, for any numbers $\alpha_2,\ldots,\alpha_m \in \bbR$ with
$|\alpha_1|,\ldots,|\alpha_m|\le \alpha_1$ there exist isometries
$\{U_i\}_{i=1}^m$ in $H$ such that $(U_i(u_i))_1= \alpha_i$
for $i=1,\ldots,m$.
Hence by Proposition~\ref{fiso} the element $v$ with
\begin{equation*}
v_i = \begin{cases}
U_i(u_i)&\text{if $i\le m$}\\
0&\text{if $i>m$}
\end{cases}
\end{equation*}
is Flinn in $X_n(H)$.
By Theorem~\ref{flinn} and Lemma~\ref{5.2} $\bar v=  
(v_{i,1})_{i=1}^{2^n}
\in {\calF}(X_n)$.
And, since the sequence $\{\alpha_2,\ldots,\alpha_m\}$ is arbitrary,  
that
implies that every element with support of cardinality smaller or  
equal
than $m$ is Flinn in $X_n$.
But if $m\ge2$ Theorem~\ref{cf} implies that $X_n = \ell_p^{2^n}$ for  
some $1\le p\le
\infty$ contrary to our assumption.
So $m=1$.
\end{pf}

As mentioned above, in the case when $H= \ell_2,$  
Proposition~\ref{5.2a} is valid for any r.i. space $X.$
Namely we have:

\begin{prop} \lb{5.5l2}
Let $X_n$ be a $n$--dimensional r.i. space not isometric to
$\ell_2^{n} \ \ (n\ge 2).$ If $u = (u_i)_{i=1}^{n} \in  
\calF(X_n(L_2)) $ then
there exists $1\le i_0 \le n$ such that $\|u_i\|_2 =0$ for all $i
\neq i_0.$
\end{prop}

\begin{rem}
We use here notation $L_2$ for the separable Hilbert space to stress  
the fact that it is nonatomic. Clearly $L_2$ is isometric to $\ell_2$  
and $X_n(L_2)$ is isometric to $X_n(\ell_2)$ via a surjective  
isometry which preserves disjointness in a vector sense and hence our  
result is valid also in $X_n(\ell_2).$
\end{rem}

\begin{pf}
Let $u \in \calF(X_n(L_2)) $ be such that $m = \card\{i : u_i
\not\equiv 0 \}$ is maximal. By Proposition~\ref{fiso} we can assume
\buo that  $u_i \equiv 0 $ for $i = m+1,\dots,n$ and $\supp u_i =
[0,1]$  for $i = 1,\dots,m.$

If we consider $X_n(L_2)$ as a function space on $\{1,\dots,n\}
\times
[0,1],$ then $\supp u_i = \{1,\dots,m\} \times  [0,1] = B.$ Since
$X_n(L_2)$ is nonatomic we can apply Theorem~\ref{kw} to conclude
that there exists a measurable function $w$ such that $\supp w = B$
and for every $x \in X_n(L_2)(B)$,
\begin{equation}\lb{n}
\|x\| = \left(\dint |x|^2w\  d\mu\right)^{1/2}\ .
\end{equation}
Since $X_n$ and $L_2$ are r.i. $w$ is constant,  say $w \equiv k.$

We need to show that $m=1.$

First notice that $m<n$ since $X_n$ is not isometric to $\ell_2^{n}$
and \eqref{n}.
Assume, for contradiction, that $m \ge 2$
and consider any element $z = (z_i)_{i=1}^{n} \in X_n(L_2)$ such
that $z_i \equiv 0$ for $i = m+2,\dots,n.$ Define $v, x, y \in
X_n(L_2) $ by

\begin{equation*}
v_i = \begin{cases}
         0 &\text{if } i\ne m+1, \\
   z_{m+1} &\text{if } i= m+1
    \end{cases};
x_i = \begin{cases}
         z_i &\text{if } i\le m, \\
         0 &\text{if } i> m
       \end{cases};
y_i = \begin{cases}
         \|x\|_2 &\text{if } i=1, \\
         0 &\text{if } i>1
       \end{cases};
\text{ resp.}
\end{equation*}

Then $\supp v \cap B = \emp, \ x, y \in X_n(L_2)(B)$ and $ \|x\|
=\|y\|$ so by Theorem~\ref{kw}$(b)$  $\|v+x\|= \|v+y\|$ i.e. $\|z\|=
\|v+y\|.$ Since $X_n$ is r.i.
$$\|v+y\| = k (\|z_{m+1}\|_2^2 + \|x\|_2^2)^{1/2} = k \|z\|_2.$$

Hence $\|z\| = k \|z\|_2$ for every $z \in X_n(L_2)(\{1,\dots,m+1\}
\times [0,1])$ and Theorem~\ref{kw} quickly leads to contradiction
with maximality of $m.$
\end{pf}

\begin{prop} \lb{disj}
Suppose that H is a sparable Hilbert space and  $X$ is a 
rearrangement-invariant   function space with
norm not proportional to the $L_2$-norm.
Suppose further that either
 $X$ is nonatomic on $[0,1]$ or it is a sequence space
$(\dim X\le \infty),$  and
\begin{enumerate}
\item[a)] $H=\ell_2;$ \/ or
\item[b)] $H=\ell_2^d$ , $X$ has a norm not proportional
to $L_p-$norm for any $1\le p\le\infty,$ $X$ satisfies property  
$(P')$  and $\dim X \ge 3$.
\end{enumerate}
Then every surjective
isometry $T:X(H) \lra X(H)$ preserves disjointness in a  vector  
sense.
\end{prop}

\begin{pf}
We will present the proof in the case when $X$ is nonatomic. If $X$
is a sequence space the proof is almost identical and slightly
simpler.

Let us denote
$e_{i,j}^n = e_i^n \otimes e_j\in X_n(H)$ ($e_j$ denotes elements of  
natural
basis of $H$) and $f_{i,j}^n = Te_{i,j}^n$ for $j, n \in \bbN,\
i \le 2^n.$

Define for any  $\om \in [0,1]
\times \bbN$
(or $\omega \in [0,1] \times \{1,\ldots,d\}$ in case (b))
$$F_n(\om) = \sum_{i=1}^{2^n} \sum_{j=1}^{\infty}
f_{i,j}^n(\om)e_{i,j}^n.$$
Following the argument same as in Theorem~6.1 of \cite{KR} we see  
that for almost every  $\omega$ \/ $F_n(\om)\in \calF(X_n'(H)).$

For the sake of completness we present this argument here.

Denote by $\Pi(X(H))$ the set of pairs $(x,x^*)$ where 
$x\in X(H), \ x^*\in X'(H)$  and $1= \|x\|=\|x^*\|=x^*(x).$

We note first that by Proposition~2.5 of \cite{KR}, $T^{-1}$
is $\sigma(X(H),X'(H))-$continuous and so has an adjoint
$S=(T^{-1})':X'(H)\lra X'(H).$ We define  
$g^n_i=Se^n_i$.  Suppose $(x,x^*)\in \Pi(X_n(H))$ where $  x=\sum
a_{i,j}e^n_{i,j}$ and $  x^*=\sum a_{i,j}^*e_{i,j}^n.$ Then  
$(Tx,Sx^*)\in
\Pi(X(H))$ and this implies that
\begin{equation}
(\sum_{i=1}^{2^n}\sum_{j=1}^{\infty}a_{i,j}f^n_{i,j}(\omega))
(\sum_{i=1}^n\sum_{j=1}^{\infty}a_{i,j}^*g^n_{i,j}(\omega))\ge
0\lb{(*)}
\end{equation}
for $\mu-$a.e.  $\omega\in\Omega.$

Using the fact that $\Pi(X_n(H))$ is separable it follows that
there is a set of measure zero $\Omega^n_0$ so that if
$\omega\notin \Omega^n_0,$ \eqref{(*)} holds for every 
$(x,x^*)\in
\Pi(X_n(H)).$ 
Let $  \Omega_0=\cup_{n\ge 1}\Omega_0^n.$

Now define
$  
G_n(\omega)=\sum_{i=1}^{2^n}\sum_{j=1}^{\infty}g^n_{i,j}(\omega)e^n_{i 
,j}\in X_n(H).$ The
above remarks show that 
if
$\omega\notin\Omega_0$  then 

$x^*(G_n(\omega))\cdot F_n(\omega)(x)\ge 0$ for all  
$(x,x^*)\in
\Pi(X_n'(H)),$ i.e.  $F_n(\omega)\in \calF(X'_n(H))$ provided that 
$G_n(\omega)\neq 0$ and $\omega\notin\Omega_0.$
We will show that this happens for a.e. $\omega\in [0,1].$

Let $B_n=\{\omega:G_n(\omega)=0\}.$ Clearly $(B_n)$ is a
descending sequence of Borel sets.  Let $B=\cap B_n.$ If
$\mu(B)>0$ then there exists a nonzero $h\in X(H)$ supported on
$B$ and $\langle h, Sx'\rangle=0$ for every $x'\in X'(H)$.
Thus $T^{-1}h=0,$ which contradicts the fact that $T$ is an  
isometry.

Let $D_n=\Omega\setminus(\Omega_0\cup B_n).$ Then $\mu(D_n)=0$ and if  

$\omega\in
D_n$ then $G_n(\omega)\neq 0$ and so it follows that
$F_n(\omega)\in\cal F(X'_n).$

Hence, by Proposition~\ref{5.2a},
\begin{equation} \lb{disj1}
\text{for a.e. $\om$ \ \ \ }
\exists i_{\om} \ \ \text{ so that   } f_{i,j}^n(\om) =0 \ \ \
\forall i \ne i_{\om}, j \in \bbN.
\end{equation}

 Let $\n_1, \n_2$ be any natural numbers. Consider the isometry $V$
of $H$ defined by
\begin{equation*}
V(e_j) = \begin{cases}
           e_j &\text{ if  } j \ne \n_1,\n_2,\\
      \frac{1}{\sqrt{2}}(e_{\n_1} + e_{\n_2}) &\text{ if  }  
j=\n_1,\\
      \frac{1}{\sqrt{2}}(e_{\n_1} - e_{\n_2}) &\text{ if  } j=\n_2,
          \end{cases}
\end{equation*}
and the induced isometry $\overline{V}$ of $X(H)$ defined by $V$
on each fiber.
\begin{equation*}
\overline{V}Te_{i,j}^n(t,\n) = \begin{cases}
           f_{i,j}^n(t,\n) &\text{ if  } \nu \ne \n_1,\n_2,\\
      \frac{1}{\sqrt{2}}(f_{i,j}^n(t,\n_1) + f_{i,j}^n(t,\n_2))
&\text{ if  } \nu =\n_1,\\
      \frac{1}{\sqrt{2}}(f_{i,j}^n(t,\n_1) - f_{i,j}^n(t,\n_2))
&\text{ if  } \nu =\n_2,
          \end{cases}
\end{equation*}

Similarly as in \eqref{disj1} we conclude that for almost every $t$  
there exists
$\bar{\imath}_{(t,v_1)}$
such that $\overline{V}Te_{i,j}^n(t,\n_1) = 0$ for all
$i \ne \bar{\imath}_{(t,v_1)}$.
Therefore, for a.e. $t,$
\begin{equation*}
f_{i,j}^n(t,\n_1) + f_{i,j}^n(t,\n_2) = 0 \ \ \forall i \ne
\bar{\imath}_{(t,v_1)}\ , \forall\ j\ .
\end{equation*}

Combining this with \eqref{disj1} we get that for almost every $t  
\in
[0,1]$ and any $\n_1,\n_2 \in \bbN\ \ $   
$\bar{\imath}_{(t,v_1)}=i_{t,\n_1} =
i_{t,\n_2}.$ It follows easily that $T$ preserves disjointness of
functions supported in disjoint dyadic intervals.
\end{pf}

We are now ready to present the main result of this paper.

\begin{th} \lb{isoxl2}

Suppose that  $X$ is a
rearrangement-in\-va\-riant  function space with norm not  
proportional to the
$L_2$-norm.
Suppose further that either
 $X$ is nonatomic on $[0,1]$ or it is a sequence space
$(\dim X\le \infty),$  and let $H$ be a separable Hilbert space.

Suppose that   $T:X(H) \lra X(H)$ is a surjective isometry. Then
there exists a nonvanishing Borel function $a $ on $\Om$ (where $\Om
= [0,1]$ if $X$ is nonatomic or $\Om \subset \bbN$ if $X$ is a
sequence space) and an invertible Borel map $\sigma:\Om\lra\Om$ such
that, for any Borel set
$B \subset \Om,$ we have
$\m(\sigma^{-1}B)=0$ if and only if $\m(B)=0$ and
a strongly measurable map $S$ of $\Om$ into $\calB(H)$
so that $S(t)$ is an isometry of $H$ onto itself for almost all $t$
and
$$Tf(t)=a(t)S(t)(f(\sigma(t))) \text{  a.e.}$$ for any $f\in X(H).$

Moreover if $X$ is not equal to $L_p[0,1]$ up to equivalent  
renorming
then $|a|= 1$ a.e. and $\s$ is measure-preserving.
\end{th}

\begin{pf}
We prove the theorem under the assumption that either:
\begin{enumerate}
\item[a)] $H=\ell_2;$ \/ or
\item[b)] $H=\ell_2^d$ , $X$ has a norm not proportional
to $L_p-$norm for any $1\le p\le\infty,$ $X$ satisfies property  
$(P')$  and $\dim X \ge 3$.
\end{enumerate}

If $\dim X = 2$ the theorem follows from Theorem~3.12 of \cite{R86}.
If $X= L_p[0,1]$, $p\ne2$ the theorem was proved by Greim \cite{G}  
and
Cambern \cite{C81}.
If $X$ does not satisfy property $(P')$ then $X'$ does and the  
result
follows by duality argument. That is, Proposition~2.5 of \cite{KR}  
says
that the isometry $T$ has an adjoint $T':X'(H)\to X'(H)$ which is a  
surjective isometry and thus has 
a canonical vector form. By Lemma~\ref{2.4'} $T''$ and hence $T$  has  
a canonical vector form.

So in the following we assume that the assertion of  
Proposition~\ref{disj}
holds, i.e., the isometry $T$ preserves disjointness.

We follow almost exactly the argument of Sourour \cite[Theorems 3.1
and 3.2]{S78}.

Let $\{x_n\}$ be the countable linearly independent subset of $H$
whose linear span $\calD$ is dense in $H$ and let $\calD_0$
be the set of all linear combinations of $\{x_n\}$ with rational
coefficients. For any measurable set $E$ let
$$\Phi(E) = \bigcup_n \supp(T(\h_E x_n)).$$

Then, since $T$ is 1-1, $\Phi$ is a set-isomorphism.

Let $y_n = T(\underline{x}_n).$ For every $t \in \Om$ define
$$A(t)x_n = y_n(t)$$
and extend $A(t)$ linearly to $\calD$ and thus for every $y \in
\calD$ $\ \ A(\cdot)y = T(\underline{y})$ a.e.
We will now extend $A(t)$ to a bounded operator on $X.$ Let  
$E\subset
\Om$ be measurable and $y \in \calD_0, $ then
\begin{equation} \lb{isoxl21}
\begin{split}
\|A(t)y \h_{\Phi(E)} \|_{X(H)} &= \| T(\underline{y}(t)\h_{\Phi(E)}
\|_{X(H)} \\
&= \| T(\underline{y}\h_{E} )\|_{X(H)} =
\|\underline{y}\h_{E}\|_{X(H)} \\
&= \|\h_E\|_X \|y\|_2
\end{split}
\end{equation}

By absolute continuity we can define for almost every $t:$
\begin{equation*}
 a(t) = \lim\begin{Sb} \m(E) \to 0 \\ t\in E \end{Sb}
\frac{\|\h_{\Phi^{-1}(E)}
 \|_X}{\|\h_E\|_X }
 \end{equation*}
(notice that if $X=L_p$ then $a(\cdot)$ coincides with the function
$h(\cdot)$ considered by Sourour).

By \eqref{isoxl21} $A(t) = a(t)S(t)$ a.e. where $S(t)$ is an  
isometry
of $H.$

The strong measurability of $S$ and surjectivity of almost all  
$S(t)$
follow as in the proof of Sourour without change.

The final remark is now an immediate consequence of
Theorem~{7.2} of \cite{KR}.
\end{pf}

\begin{ack}
I wish to express my gratitude to Professor Nigel Kalton for his  
interest
in this work and many valuable disscussions.
\end{ack}
 \bibliographystyle{standard}
    \bibliography{tref}

\end{document}